\documentclass[10pt]{article}
\usepackage{latexsym}
\usepackage{amsfonts,amsmath,amssymb}

\usepackage[english,frenchb]{babel}
\usepackage[latin1]{inputenc}
\usepackage{t1enc}

\newtheorem{defi} {D{\'e}finition}
\newtheorem{theo} {Th{\'e}or{\`e}me}   
\newtheorem{prop} {Proposition}

\newtheorem{lem} {Lemme}
\newtheorem{Rmq} {Remarque}

\begin{document}

\selectlanguage{frenchb}

     \bibliographystyle{smfplain}

  \begin{center}  
    
    \Large  \bf

K-théorie {\'e}quivariante des vari{\'e}t{\'e}s de drapeaux et des variétés de
Bott-Samelson.

\end{center}

\rm
\normalsize

 \vspace {0,25 cm}

 \begin{center}

\large

Matthieu WILLEMS

\end{center}

\rm
\normalsize

\vspace{0,3 cm}

\begin{abstract}

Le but de ce texte est de donner explicitement les valeurs des
restrictions aux points fixes d'une base $\hat{\psi}_{w}$
(resp. $\hat{\mu}_{\epsilon}$) de la $K$-th{\'e}orie
{\'e}quivariante des vari{\'e}t{\'e}s de drapeaux (resp. des variétés de
Bott-Samelson). On donne tout d'abord une
démonstration combinatoire de la formule des $\psi_{w}$ (théorème $1$).
Grace à la formule de localisation, on calcule ensuite une base $\mu_{\epsilon}$
de la $K$-théorie
équivariante des variétés de Bott-Samelson (théorème $3$) ce qui nous donne une
démonstration plus géométrique du théorème $1$. 
Cette étude nous permet également de calculer la matrice de
changement de bases entre $\hat{\psi}_{w}$ et $*[\mathcal{O}_{\overline{X}_{w}}]$.

\end{abstract}

\selectlanguage{english}

\begin{abstract}

The aim of this text is to give an explicit formula  for the restriction to the
fixed points of a basis  $\hat{\psi}_{w}$
(resp. $\hat{\mu}_{\epsilon}$) of the equivariant $K$-theory of the flag
varieties (resp. of the Bott-Samelson varieties). First of all, we give a
combinational proof of the formula for the $\psi_{w}$ (theorem $1$). Then, we calculate 
a basis $\mu_{\epsilon}$ of the equivariant $K$-theory of the Bott-Samelson 
varieties using the localization formula (theorem $3$). Then we give a more
geometric proof of the theorem $1$ using the theorem $3$. In the finite case, we
describe how the basis $*[\mathcal{O}_{\overline{X}_{w}}]$ transforms with respect to the
basis $\hat{\psi}_{w}$.

\end{abstract}

\selectlanguage{frenchb}

\vspace{1 cm}

Après avoir rédigé cette note, j'ai eu connaissance de résultats de William
Graham qui prouve la formule du théorème $1$ de deux manières différentes
dans le preprint \cite{gra}. Une de ses deux démonstrations utilise des idées
similaires à celles développées dans les sections $3$ et $4$. Ces idées sont
utilisées par Sarah Billey dans \cite{sb} dans le cas de la
cohomologie équivariante. Pour le type $A$, les algèbres de Hecke sont également
utilisées dans les articles \cite{llt} et \cite{fk2} dans le cadre de la 
$K$-théorie équivariante. 
Dans \cite{gra}, William Graham donne des formules explicites pour la
restriction aux points fixes des classes $[\mathcal{O}_{\overline{X}_{w}}]$.

\section{Préliminaires et notations}

Les définitions et les résultats qui suivent sur les algèbres de Kac-Moody
 sont exposés dans \cite{km}. Soit
 $A=(a_{ij})_{1\leq i,j \leq r}$ une matrice de Cartan généralisée (c'est-à-dire 
telle que $a_{ii}=2$, $-a_{ij} \in \mathbb{N}$ si $i \neq j$, et
  $a_{ij}=0$ si et seulement si $a_{ji}=0$). On choisit un triplet $(\mathfrak{h},
  \mathfrak{\pi}, \mathfrak{\pi^{\vee}})$ (unique à isomorphisme près), où
  $\mathfrak{h}$ est un $\mathbb{C}$-espace vectoriel de dimension $(2r-{ \rm
 rg}(A))$,
  $\mathfrak{\pi} = \{\alpha_{i}\}_{1 \leq i \leq r} \subset \mathfrak{h}^*$, et 
$\mathfrak{\pi^{\vee}} = \{h_{i}\}_{1 \leq i \leq r} \subset \mathfrak{h}$ sont
  des ensembles d'éléments linéairement indépendants vérifiant
  $\alpha_{j}(h_{i})=a_{ij}$. L'algèbre de Kac-Moody
  $\mathfrak{g}=\mathfrak{g}(A)$ est l'algèbre de Lie sur $\mathbb{C}$ engendrée
  par $\mathfrak{h}$ et par les symboles $e_{i}$ et $f_{i}$ ($1 \leq i \leq r$)
  soumis aux relations $[\mathfrak{h},\mathfrak{h}]=0$, $[h,
  e_{i}]=\alpha_{i}(h)e_{i}$, $[h, f_{i}]=-\alpha_{i}(h)f_{i}$ pour tout $h \in
  \mathfrak{h}$ et tout $1 \leq i \leq r$, $[e_{i}, f_{j}]=\delta_{ij}h_{j}$
  pour tout $1 \leq i,j \leq r$, et : 
$$({\rm ad }e_{i})^{1-a_{ij}}(e_{j})=0=({\rm ad }f_{i})^{1-a_{ij}}(f_{j})
   \hspace{0,2 cm} \forall \hspace{0,2 cm}  1 \leq i \neq j \leq r.$$

L'algèbre $\mathfrak{h}$ s'injecte canoniquement dans $\mathfrak{g}$. On
l'appelle la sous-algèbre de Cartan de $\mathfrak{g}$. On a la
décomposition suivante : 
$$\mathfrak{g}=\mathfrak{h} \oplus \sum_{\alpha \in
  \Delta_{+}}(\mathfrak{g}_{\alpha} \oplus \mathfrak{g}_{-\alpha}),$$
où pour $\lambda \in \mathfrak{h}^*$, $\mathfrak{g}_{\lambda} = \{ x \in
\mathfrak{g} \: {\rm tels \: que }\: [h, x]=\lambda(h)x, \forall h \in
\mathfrak{h} \}$, et où on définit $\Delta_{+}$ par $\Delta_{+} = \{ \alpha \in
\sum_{i=1}^{r}\mathbb{N}\alpha_{i} \: {\rm tels \: que } \: \alpha \neq 0
\:{\rm et } \: \mathfrak{g}_{\alpha} \neq 0 \}$. On pose $\Delta=\Delta_{+} \cup
\Delta_{-}$ où $\Delta_{-} = -\Delta_{+}$. On appelle $\Delta_{+}$
(respectivement $\Delta_{-}$) l'ensemble des racines positives (respectivement
négatives). Les racines $\{\alpha_{i}\}_{1 \leq i \leq r}$ sont appelées les
racines simples. On définit une sous-algèbre de Borel $\mathfrak{b}$ de $\mathfrak{g}$ par
$\mathfrak{b}=\mathfrak{h} \oplus \sum_{\alpha \in \Delta_{+}}
\mathfrak{g}_{\alpha}$.

Au couple $(\mathfrak{g}, \mathfrak{h})$, on associe le groupe de Weyl $W
\subset { \rm Aut}(\mathfrak{h}^*)$, engendré par les réflexions simples $\{r_{i}\}_{1
  \leq i \leq r}$, où $r_{i}(\lambda)=\lambda-\lambda(h_{i})\alpha_{i}$ pour
tout $\lambda \in \mathfrak{h}^*$. Le groupe $W$ étant un groupe de Coxeter, on
a une notion d'ordre de Bruhat qu'on notera $u \leq v$ et une notion de longueur
qu'on notera $l(w)$. On notera $1$ l'élément neutre de $W$ et dans le cas fini
(i.e $W$ fini $\Leftrightarrow$ $A$ définie positive $\Leftrightarrow
\mathfrak{g}$ de dimension finie), 
on note $w_{0}$ le plus grand élément de $W$. 
Le groupe de Weyl préserve $\Delta$. On pose $R=W\pi$ et $R^{+}=R\cap
\Delta_{+}$. Pour $\beta = w\alpha_{i} \in R^{+}$, on pose
$r_{\beta}=wr_{i}w^{-1} \in W$ (qui est indépendant du choix du couple $(w,
\alpha_{i})$ vérifiant $\beta = w\alpha_{i}$).
Pour un élément $w$ de
$W$, on définit l'ensemble $\Delta(w)$ des inversions de $w$ par
 $\Delta(w)=\Delta_{+} \cap w^{-1}\Delta_{-}$.

On fixe un réseau $\mathfrak{h}_{\mathbb{Z}} \subset \mathfrak{h}$ tel que : 

\vspace{0,1 cm}

(i) $\mathfrak{h}_{\mathbb{Z}} \otimes_{\mathbb{Z}}\mathbb{C}=\mathfrak{h}$,

(ii) $h_{i} \in \mathfrak{h}_{\mathbb{Z}}$ pour tout $1 \leq i \leq r$,

(iii) $\mathfrak{h}_{\mathbb{Z}}/ \sum_{i=1}^{r}\mathbb{Z}h_{i}$ est sans
torsion,

(iv) $\alpha_{i} \in \mathfrak{h}_{\mathbb{Z}}^* =
{ \rm Hom}(\mathfrak{h}_{\mathbb{Z}}, \mathbb{Z})$ ($\subset \mathfrak{h}^*$) pour tout 
$1 \leq i \leq r$. 

\vspace{0,2 cm}

On choisit des poids fondamentaux $\rho_{i} \in \mathfrak{h}_{\mathbb{Z}}^*$ 
($1 \leq i \leq r$) qui vérifient $\rho_{i}(h_{j})=\delta_{ij}$, pour tout 
$1 \leq i,j \leq r$. On pose $\rho=\sum_{i=1}^{r}\rho_{i}$.

On note $G=G(A)$ le groupe de Kac-Moody associé à $\mathfrak{g}$ par Kac et
Peterson dans \cite{kp}. Dans le cas fini, $G$ est un groupe de Lie semi-simple
complexe connexe et simplement connexe. On note $H \subset B \subset G$ les sous-groupes de $G$
associés respectivement à $\mathfrak{h}$ et $\mathfrak{b}$. Soit $K$ la forme
unitaire standard de $G$ et $T=K \cap H$ le tore maximal de $K$ associé à
$\mathfrak{h}$.  On pose $X=G/B=K/T$. On fait agir $T$ sur $X$ par
multiplication à gauche.

Soit $X[T]$ le groupe des caractères de $T$, on pose
$R[T]=\mathbb{Z}[X[T]]$ et on note $Q[T]$ le corps des fractions de $R[T]$. Pour
un poids entier $\lambda$, on note $e^{\lambda} \in X[T]$ le caractère
correspondant. 

On note $F(W,R[T])$ (respectivement $F(W,Q[T])$) l'algèbre des fonctions sur $W$ à
valeurs dans $R[T]$ (respectivement $Q[T]$) munie de
l'addition et de la multiplication point par point. Pour tout $1 \leq i \leq r$,
on définit un opérateur de Demazure $D_{i}$ sur $F(W,Q[T])$ par : 

$$(D_{i}f)(v)=\frac{f(v)-f(vr_{i})e^{-v\alpha_{i}}}{1-e^{-v
\alpha_{i}}}.$$
   
Les opérateurs de Demazure vérifiant les relations de tresses de $W$,
on peut définir un opérateur $D_{w}$ pour tout $w \in W$. On note $\Psi$ la
sous-algèbre de $F(W,R[T])$ définie par : 
$$ \Psi = \{ f \in F(W,R[T]), \hspace{0,2 cm} {\rm telles \hspace{0,2 cm} que } \hspace{0,2 cm}
\forall w \in W, \hspace{0,2 cm} D_{w}f \in F(W,R[T]) \}.$$

\section{$K$-théorie équivariante des variétés de drapeaux} 

On définit la $K$-théorie $T$-équivariante de $X=G/B$ comme le groupe construit à
partir du semi-groupe des classes
d'isomorphisme de fibrés complexes $T$-équivariants au dessus de $X$. On munit ce groupe
d'une structure d'anneau définie à l'aide du produit tensoriel. De
plus comme la  $K$-théorie $T$-équivariante du point s'identifie à $R[T]$, on
obtient une structure de 
$R[T]$-algèbre qu'on notera $K_{T}(X)$. L'ensemble $X^T$ des points fixes de $X$ 
sous l'action de $T$ s'identifie à
$W$. L'injection de $X^T$ dans $X$ définit une application $i_{T}^* : K_{T}(X)
\rightarrow  K_{T}(X^T)$. De plus, l'ensemble des points fixes étant discret, on
peut identifier $K_{T}(X^T)$ avec $F(W,R[T])$ et on a ainsi une application 
$i_{T}^* : K_{T}(X) \rightarrow F(W,R[T])$. On notera $*$ l'involution de
$K_{T}(X)$ définie par la dualité des fibrés et on notera de la m\^eme façon
l'involution de $R[T]$ définie sur les caractères par
$*(e^{\lambda})=e^{-\lambda}$, ce qui induit une involution de $F(W,R[T])$.
 Pour tout élément $\tau \in K_{T}(X)$, 
$*i_{T}^*(\tau)=i_{T}^*(*\tau)$. Le r{\'e}sultat suivant est prouv{\'e} dans \cite{kkk} :

\begin{prop}

L'application $i_{T}^*$ est injective et l'image de $K_{T}(X)$ par cette
application est égale à $\Psi$. De plus, $\Psi=\prod_{w \in W}R[T]\psi^{w}$, où
les fonctions $\psi^{w}$ sont caractérisées par les propriétés suivantes :  

 (i) $\psi^{w}(v)=0$ sauf si $w \leq v$,
   
   (ii) $\psi^{w}(w) = \prod_{\beta \in \Delta(w^{-1})}
(1-e^{\beta})$,

   $ (iii)\left\{ \begin{array}{ll} D_{i}\psi^{w} = \psi^{w} + \psi^{wr_{i}}
 & {\rm si } \hspace{0.2 cm} wr_{i}<w, \\ 
D_{i}\psi^{w}=0
 & {\rm si } \hspace{0.2 cm} wr_{i}>w,  \end{array} \right.$

(iv) $ \forall v \in W, \psi^{1}(v) =e^{\rho-v\rho} $.
    
\end{prop}

On pose $\hat{\psi}^w=(i_{T}^*)^{-1}(\psi^{w})$.

\begin{Rmq}

Un élément $f=(a_{w})_{w \in W}$ de $\prod_{w \in W}R[T]\psi^{w}$ est bien une 
fonction de $W$ à valeurs dans $R[T]$. En effet soit $v \in W$, d'après la propriété
$(i)$, $\sum_{w \in W}a_{w}\psi^{w}(v)$ est une somme finie où les termes
éventuellement non nuls correspondent aux éléments $u$ de $W$ qui vérifient $u
\leq v$.

\end{Rmq}

 Dans \cite{kkk}, B. Kostant et S. Kumar composent $i_{T}^{*}$ avec $\phi : F(W, Q[T]) \rightarrow
F(W, Q[T])$ définie par 
$\phi(f)(w)=f(w^{-1})$ pour tout élement $f$ de $F(W, Q[T])$ et tout $w \in
W$. Ils trouvent alors la sous algèbre $\Psi'$ (notée $\Psi$ dans \cite{kkk}) 
de $F(W,R[T])$ : 
$$\Psi' = \{ f \in F(W,R[T]), \hspace{0,2 cm} {\rm telles \hspace{0,2 cm} que } \hspace{0,2 cm}
\forall w \in W, \hspace{0,2 cm} D'_{w}f \in F(W,R[T]) \},$$
où les opérateurs $D_{w}'$ sont définis à partir des opérateurs $D_{i}'$ donnés
par : 
$$(D_{i}'f)(v)=\frac{f(v)-f(r_{i}v)e^{-v^{-1}\alpha_{i}}}{1-e^{-v^{-1}
\alpha_{i}}}.$$

Ils considèrent la base $\psi_{w} '$ (notée $\psi^{w}$ dans \cite{kkk}) 
de $\Psi '$ reliée à la base $\psi^{w}$ de la proposition $1$ par la relation 
$\psi_{w}'=\phi(\psi^{w^{-1}})$. Pour tout couple $(w,v) \in W^2$, $\psi_{w}'(v)=
\psi^{w^{-1}}(v^{-1})$.

On définit le mono\"{i}de $\underline{W}$ comme le mono\"{i}de engendré par les éléments
$\{\underline{r}_{i}\}_{1 \leq i \leq r}$ soumis aux relations
$\underline{r}_{i}^2=\underline{r}_{i}$ et
aux m\^emes relations de tresses que les éléments $r_{i}$ de $W$. D'après l'étude
générale des algèbres de Hecke, l'ensemble $\underline{W}$ s'identifie à l'ensemble
$W$. Pour un élément $w$ de $W$, on notera $\underline{w}$ l'élément correspondant
dans $\underline{W}$ et pour $\underline{v} \in \underline{W}$, on notera $v$ l'élément
associé dans $W$. Dans $\underline{W}$, on a les relations suivantes : 

\begin{equation} \label{hecke1}  \left\{ \begin{array}{ll}
      \underline{w}\hspace{0,1 cm}\underline{r}_{i}=
\underline{wr_{i}}
 & {\rm si } \hspace{0.2 cm} wr_{i}>w, \\ 
\underline{w} \hspace{0,1 cm}\underline{r}_{i}=\underline{w}
 & {\rm si } \hspace{0.2 cm} wr_{i}<w.  \end{array} \right.
\end{equation}

\begin{equation} \label{hecke2}
 \left\{ \begin{array}{ll} \underline{r}_{i}\hspace{0,1 cm}\underline{w}=\underline{r_{i}w}
 & {\rm si } \hspace{0.2 cm} r_{i}w>w, \\ 
\underline{r}_{i} \hspace{0,1 cm} \underline{w}=\underline{w}
 & {\rm si } \hspace{0.2 cm} r_{i}w<w.  \end{array} \right.
\end{equation}

Soit $v \in W$ et soit $v=r_{i_{1}} \cdots r_{i_{l}}$ une décomposition réduite de
$v$. Pour $1\leq j \leq l$, on d{\'e}finit un 
{\'e}l{\'e}ment $\beta_{j} \in 
\mathfrak{h}^*$ par $\beta_{j}=r_{i_{1}} \cdots r_{i_{j-1}}
\alpha_{i_{j}}$ .

\begin{theo}

 Si $w \in W$ est tel que $w \leq v$, on a la formule suivante :    

   $$ \psi^{w}(v) =e^{\rho - v\rho}\sum_{l(w) \leq m \leq l(v)}
 \sum
 (e^{-\beta_{j_{1}}}-1) \cdots (e^{-\beta_{j_{m}}}-1),
$$
    o{\`u} la deuxi{\`e}me somme porte sur l'ensemble des entiers $1\leq j_{1} 
   < \cdots < j_{m} \leq l$ tels que $\underline{r_{i_{j_{1}}}} \ldots   
  \underline{r_{i_{j_{m}}}}=\underline{w}$.
    
\end{theo}

Donnons quelques exemples de calculs pour expliciter cette 
formule.

\vspace {0,1 cm}

Tout d'abord pour tout $v\in
W$, on retrouve bien $\psi^{1}(v)=e^{\rho-v\rho}$, puisque la seule façon
de trouver $1$ en dessous de $v$ est de prendre la suite vide.

\vspace {0,2 cm}

Plaçons nous dans le cas où $G=SL_{4}(\mathbb{C})$. Calculons $\psi^{w}(v)$
avec $w=r_{3}r_{2}$ et $v=r_{2}r_{3}r_{2}r_{1}r_{2}$. Il y a $3$ façons de
``trouver $\underline{w}$ en dessous de $v$'' :
$\underline{w}=\underline{r_{i_{2}}}\hspace{0,1 cm} \underline{r_{i_{3}}}$, 
$\underline{w}=\underline{r_{i_{2}}}\hspace{0,1 cm} \underline{r_{i_{5}}}$, 
$\underline{w}=\underline{r_{i_{2}}}\hspace{0,1 cm}
\underline{r_{i_{3}}}\hspace{0,1 cm}\underline{r_{i_{5}}}$, et on trouve donc :
$$\psi^{w}(v)=e^{\alpha_{2}+(\alpha_{2}+\alpha_{3})+\alpha_{3}+(\alpha_{1}+\alpha_{2}+\alpha_{3})+(\alpha_{1}+\alpha_{2})}[(e^{-(\alpha_{2}+\alpha_{3})}-1)(e^{-\alpha_{3}}-1)$$
$$+
(e^{-(\alpha_{2}+\alpha_{3})}-1)(e^{-(\alpha_{1}+\alpha_{2})}-1)  +
(e^{-(\alpha_{2}+\alpha_{3})}-1)(e^{-\alpha_{3}}-1)(e^{-(\alpha_{1}+\alpha_{2})}-1)]$$
$$=e^{2\alpha_{1}+4\alpha_{2}+3\alpha_{3}}[1+e^{-(\alpha_{1}+2\alpha_{2}+2\alpha_{3})}-
e^{-(\alpha_{2}+\alpha_{3})}-
e^{-(\alpha_{1}+\alpha_{2}+\alpha_{3})}]$$
$$=e^{2\alpha_{1}+4\alpha_{2}+3\alpha_{3}}+e^{\alpha_{1}+2\alpha_{2}+\alpha_{3}}-
e^{2\alpha_{1}+3\alpha_{2}+2\alpha_{3}}-e^{\alpha_{1}+3\alpha_{2}+2\alpha_{3}}.$$

\vspace {0,2 cm}

Nous allons maintenant montrer que l'expression du théorème $1$ est indépendante de
la décomposition réduite de $v$ choisie puis que les éléments de $F(W,R[T])$
ainsi définis vérifient les propriétés de la proposition $1$.

\section{Indépendance par rapport au choix d'une décomposition réduite de v}

Soit $A$ un anneau commutatif. On définit l'algèbre $\mathcal{H}$ comme la
$A$-algèbre de Hecke engendrée par $\{u_{i}\}_{1 \leq i \leq r}$ soumis aux relations de
tresses définissant $W$ et aux relations $u_{i}^2=u_{i}$. Soit $w \in W$, on
peut définir $u_{w} \in \mathcal{H}$ par $u_{w}=u_{i_{1}} \cdots u_{i_{l}}$ où
$w=r_{i_{1}}
\cdots r_{i_{l}}$ est une décomposition réduite quelconque de $w$. Les éléments 
$\{u_{w}\}_{w \in W}$ forment une base du $A$-module $\mathcal{H}$.

Soit $r_{i_{1}}, \ldots , r_{i_{k}}$ une suite de réflexions simples et soit
$\underline{w}=\underline{r}_{i_{1}} \cdots \underline{r}_{i_{k}} \in
\underline{W}$. 
D'après les relations vérifiées par les $u_{i}$, $u_{i_{1}} \cdots u_{i_{k}}=u_{w}$.

 Pour tout $1 \leq i \leq r$, on définit la fonction
 $$h_{i} : \begin{array}{l}A \rightarrow \mathcal{H} \\
x \mapsto 1 + (x-1)u_{i}.\end{array}$$

 On vérifie que ces fonctions $h_{i}$
satisfont les relations suivantes (énoncées sous une forme différente dans
\cite{fk}) : 

\begin{prop} 

Soient $1 \leq i,j \leq r$ des entiers distincts, deux éléments quelconques $x$
et $y$ de $A$ vérifient les équations suivantes :  

$$\left\{ \begin{array}{llll} h_{i}(x)h_{j}(y)=h_{j}(y)h_{i}(x)
 & {\rm si } \hspace{0.2 cm} (r_{i}r_{j})^2 = 1, \\ 
 h_{i}(x) h_{j}(xy)  h_{i}(y)=h_{j}(y)h_{i}(xy)h_{j}(x)
 & {\rm si } \hspace{0.2 cm}(r_{i}r_{j})^3 = 1, \\ 
 h_{i}(x) h_{j}(xy) h_{i}(xy^2) h_{j}(y)=h_{j}(y)h_{i}(xy^2) h_{j}(xy) h_{i}(x)
 & {\rm si } \hspace{0.2 cm} (r_{i}r_{j})^4 = 1, \\ 
h_{i}(x) h_{j}(x^3y) h_{i}(x^2y) h_{j}(x^3y^2)h_{i}(xy) h_{j}(y)
 &   \\ 
=h_{j}(y) h_{i}(xy) h_{j}(x^3y^2) h_{i}(x^2y) h_{j}(x^3y) h_{i}(x)
 & {\rm si } \hspace{0.2 cm}  (r_{i}r_{j})^6 = 1.

\end{array} \right.$$

\end{prop}

\vspace{0,2 cm}

Dans la suite, on prendra pour $A$ l'anneau $R[T]$.
Soit $w \in W$ et $w=r_{i_{1}}\cdots r_{i_{l}}$ une décomposition réduite de
$w$. On définit un élément de $\mathcal{H}$ par : 

$$\mathcal{R}_{i_{1},\ldots , i_{l}} = \prod_{j=1}^{l}h_{i_{j}}(e^{-\beta_{i_{j}}}).$$

A l'aide de la proposition $2$, une démonstration analogue à celle donnée par
S. Billey dans \cite{sb} dans le cas de l'algèbre nil-Coxeter nous donne le
résultat suivant :

\begin{theo}

Soit $w \in W$. L'élément $\mathcal{R}_{i_{1},\ldots , i_{l}}$ de $\mathcal{H}$ est
indépendant du choix d'une décomposition réduite $w=r_{i_{1}} \cdots r_{i_{l}}$ de
$w$. Il ne dépend que de $w$ et on le notera donc $\mathcal{R}_{w}$.

\end{theo}

Donnons une idée de la démonstration. D'après la définition de
$\mathcal{R}_{i_{1},\ldots , i_{l}}$ et d'après la connexité du graphe des
décompositions réduites de $w$, on peut se contenter de regarder ce qui se passe
pour un élément $w$ correspondant à une relation de tresses. Prenons par
exemple $w=r_{i}r_{j}r_{i}=r_{j}r_{i}r_{j}$. Alors $\mathcal{R}_{i, j ,
  i}=h_{i}(e^{-\alpha_{i}})h_{j}(e^{-\alpha_{i}-\alpha_{j}})h_{i}(e^{-\alpha_{j}})$ et 
$\mathcal{R}_{j, i ,
  j}=h_{j}(e^{-\alpha_{j}})h_{i}(e^{-\alpha_{i}-\alpha_{j}})h_{j}(e^{-\alpha_{i}})$, et en
utilisant la deuxième relation de la proposition $2$, on obtient le
résultat. Les autres cas se traitent de la m\^eme manière.

Le terme $\sum_{l(w) \leq m \leq l(v)} \sum (e^{-\beta_{j_{1}}}-1) 
\cdots (e^{-\beta_{j_{m}}}-1)$ du théorème $1$ est le coefficient de $R_{v}$ sur
$u_{w}$ dans la base $\{u_{w}\}_{w \in W}$ de $\mathcal{H}$ et est donc bien indépendant de la
décomposition réduite de $v$ choisie.

\section{Démonstration du théorème}

Notons $\tilde{\psi}^{w}$ l'élément de $F(W,R[T])$ défini par la formule du
théorème $1$ (on pose $\tilde{\psi}^{w}(v)=0$ si $v$ n'est pas plus grand que
$w$). Pour démontrer ce théorème, il suffit de montrer que 
les fonctions $(\tilde{\psi}^{w})_{w \in W}$ vérifient les quatre propriétés
de la proposition $1$. Les propriétés $(i)$ et $(iv)$ sont immédiates.

\vspace{0,2 cm}

Pour démontrer la propriété $(ii)$, rappelons tout d'abord les deux lemmes
suivants (voir \cite{bou}) : 

\begin{lem}

Soit $v \in W$ et $v=s_{i_{1}} \cdots s_{i_{k}}$ une décomposition réduite de $v$,
alors $ \Delta(v^{-1})=\{ \beta_{j}, 1 \leq j
  \leq k \}$.

\end{lem}

\begin{lem}

Soit $v \in W$ et $v=s_{i_{1}} \cdots s_{i_{k}}$ une décomposition réduite de $v$,
alors $\displaystyle{\rho - v\rho = \sum_{j=1}^{k}\beta_{j}}$.

\end{lem}

On a donc : 

\begin{equation}  \label{rho} e^{\rho - v\rho}=\prod_{j=1}^{l}e^{\beta_{j}}.  \end{equation}

De cette formule, on déduit pour tout $w \in W$ : 

$$\tilde{\psi}^{w}(w)=\prod_{\beta \in \Delta(w^{-1})}e^{\beta}\prod_{\beta
\in \Delta(w^{-1})}(e^{-\beta}-1)=\prod_{\beta \in
\Delta(w^{-1})}(1-e^{\beta})=\psi^{w}(w),$$
ce qui nous donne la propriété $(ii)$.

\vspace{0,2 cm}

Montrons maintenant que les $(\tilde{\psi}^{w})_{w \in W}$ vérifient la propriété $(iii)$ de la
proposition $1$.

Soit $w \in W$ et $r_{i}$ une réflexion simple. Supposons tout d'abord
$wr_{i}>w$. Il faut alors montrer que pour tout $v \in W$, on a : 

$$\tilde{\psi}^{w}(v)=\tilde{\psi}^{w}(vr_{i})e^{-v\alpha_{i}}.$$

On peut supposer $vr_{i}>v$. Si $v$ n'est pas plus grand que $w$, $vr_{i}$ non
plus car $w$ n'a pas de décomposition qui commence par $r_{i}$ car
$wr_{i}>w$. On suppose donc $w \leq v <vr_{i}$. Comme $w$ n'a aucune
décomposition qui finit par $r_{i}$, la somme est la m\^eme à gauche et à
droite de l'égalité. Il suffit donc de vérifier $e^{\rho - v\rho}=
e^{\rho - vr_{i}\rho}e^{-v\alpha_{i}}$, ce qui est une conséquence
immédiate de la formule \ref{rho}.

Supposons maintenant $wr_{i}<w$. Il faut montrer que pour tout $v \in 
W$, on a :

    \begin{equation} \label{dem} \frac{\tilde{\psi}^{w}(v)-\tilde{\psi}^{w}
(vr_{i})e^{-v\alpha_{i}}}
    {1-e^{-v\alpha_{i}}}=\tilde{\psi}^{w}(v)+\tilde{\psi}^{wr_{i}}(v). 
    \end{equation}

Supposons tout d'abord $vr_{i}>v$. On se place dans le cas o\`u $w \leq 
vr_{i}$ (sinon le r\'esultat est trivial). On choisit une d\'ecomposition 
r\'eduite $v=r_{i_{1}} \cdots r_{i_{l}}$ de $v$. On prend pour 
$vr_{i}$ la d\'ecomposition $vr_{i}=r_{i_{1}} \cdots r_{i_{l}}r_{i}$. 
On trouve alors (en utilisant la formule \ref{rho}) : 

$$\tilde{\psi}^{w}(vr_{i})e^{-v\alpha_{i}}=\tilde{\psi}^{w}(v)+
(e^{-v\alpha_{i}}-1)(\tilde{\psi}^{w}(v)+\tilde{\psi}^{wr_{i}}(v)),$$
le premier terme venant des sous d\'ecompositions de $v$ ``\'egales'' \`a 
$w$, le deuxi\`eme des m\^emes sous d\'ecompositions de $v$ auxquelles 
on rajoute $r_{i}$ à la fin et qui redonnent donc $w$ (car 
$wr_{i}<w$), et le troisi\`eme des sous d\'ecompositions de $v$ 
``\'egales'' \`a $wr_{i}$. On trouve alors bien la formule \ref{dem}.

Supposons maintenant $vr_{i}<v$. On peut appliquer ce qui pr\'ec\`ede 
\`a $v'=vr_{i}$ car $v'r_{i}>v'$ et on trouve : 

$$\frac{\tilde{\psi}^{w}(vr_{i})-\tilde{\psi}^{w}(v)e^{v\alpha_{i}}}
     {1-e^{v\alpha_{i}}}=\tilde{\psi}^{w}(vr_{i})+
     \tilde{\psi}^{r_{i}w}(vr_{i}). $$

De plus, on peut appliquer le cas $wr_{i}>w$ \`a $w'=wr_{i}$ et on 
obtient : $\tilde{\psi}^{wr_{i}}(vr_{i}) = e^{v\alpha_{i}}
\tilde{\psi}^{wr_{i}}(v)$. En substituant ainsi 
$\tilde{\psi}^{wr_{i}}(vr_{i})$ dans l'expression pr\'ec\'edente, on 
obtient la formule \ref{dem}.

\section{$K$-théorie équivariante des variétés de Bott-Samelson}

Dans toute la suite, on se place dans le cas fini; $G$ est donc un groupe de Lie
semi-simple complexe simplement connexe d'algèbre de Lie $\mathfrak{g}$, $K$
est une forme réelle compacte de $G$, $H$ est un sous-groupe de Cartan de $G$ et
$T=K \cap H$ est un tore maximal de $K$. On note $e$ l'élément neutre de $K$.
Soit $N$ un entier strictement positif. Consid{\'e}rons une suite de $N$ racines
simples $\mu_{1}$, \ldots , $\mu_{N}$ non n{\'e}cessairement distinctes. Pour 
$1 \leq i \leq N$, on note
$G_{i}$ le sous-groupe ferm{\'e} connexe de $G$ d'alg{\`e}bre de Lie
$\mathfrak{g_{\mu}}_{i}\oplus \mathfrak{h}\oplus
\mathfrak{g_{-\mu}}_{i}$ et on pose  $K_{i}=G_{i}\bigcap K$ . On d{\'e}finit : 

$$\Gamma(\mu_{1}, \ldots ,\mu_{N})=K_{1} \times_{T} K_{2} \times_{T} \cdots \times_{T}
K_{N}/T,$$
comme l'espace des orbites de $K_{1} \times K_{2} \times \cdots \times K_{N}$
sous l'action à droite de $T^N$ d{\'e}finie par : 

$$(k_{1}, k_{2}, \ldots , k_{N})
(t_{1}, t_{2}, \ldots , t_{N}) = 
(k_{1}t_{1},t_{1}^{-1} k_{2}t_{2}, \ldots ,t_{N-1}^{-1} k_{N}t_{N}),\hspace{0,1
  cm} t_{i} \in
T, \hspace{0,1 cm} k_{i} \in K_{i}.$$

On notera $[k_{1}, k_{2}, \ldots , k_{N}]$ la classe de $(k_{1}, k_{2},
\ldots , k_{N})$ dans $\Gamma(\mu_{1}, \ldots ,\mu_{N})$. On notera
$k_{\mu_{i}}$ un repr{\'e}sentant quelconque de la reflexion $r_{\mu_{i}}$ de
$N_{K_{i}}(T)/T$. Dans la suite, on
notera $\Gamma(\mu_{1}, \ldots ,\mu_{N})$ par $\Gamma$. On munit $\Gamma$
de sa structure complexe canonique d{\'e}finie dans \cite{hbs}.

On d{\'e}finit une action à gauche de $T$ sur $\Gamma$ par : 

$$t[k_{1}, \ldots ,k_{N}]=[tk_{1}, \ldots ,k_{N}],\hspace{0,1 cm} t \in T,
\hspace{0,1 cm}  k_{i} \in K_{i}.$$

On pose $\mathcal{E} =
\{0,1\}^N$. Pour $\epsilon \in \mathcal{E}$, on note $Y_{\epsilon}
\subset \Gamma$ l'ensemble des classes $[k_{1},
k_{2}, \ldots , k_{N}]$ qui v{\'e}rifient pour tout entier $i$ compris entre
$1$ et $N$ : 

$$\left\{ \begin{array}{ll} k_{i} \in T
 & { \rm si} \hspace{0,15 cm} \epsilon_{i} =0, \\ 
 k_{i} \notin T
 & { \rm si }\hspace{0,15 cm} \epsilon_{i} =1.
\end{array}\right.$$

On v{\'e}rifie imm{\'e}diatement que cette d{\'e}finition est bien compatible avec
l'action de $T^N$.  On munit $\mathcal{E}$ d'une structure de groupe en
 identifiant $\{0,1\}$ avec $\mathbb{Z}/2\mathbb{Z}$.
 Pour $\epsilon \in \mathcal{E}$, on note
$\pi_{+}(\epsilon)$ l'ensemble des entiers $i$ tels que
$\epsilon_{i}=1$ et $\pi_{-}(\epsilon)$ l'ensemble des entiers $i$
tels que $\epsilon_{i}=0$. On pose $l(\epsilon) =
{\rm card}(\pi_{+}(\epsilon))$. On note $(i) \in \mathcal{E}$
 l'{\'e}l{\'e}ment de $\mathcal{E}$ d{\'e}fini par $(i)_{j}=\delta_{i,j}$.
 Pour $\epsilon \in \mathcal{E}$, on pose
$\displaystyle{v_{i}(\epsilon) =\prod_{ {\tiny \begin{array}{ll}  1\leq k \leq i, \\ 
 k \in \pi_{ +}(\epsilon) \end{array}} 
  }r_{\mu_{k}}}$, ($v_{i}(\epsilon)=1$, si $\{1\leq k \leq i,  k \in \pi_{
+}(\epsilon)\}= \emptyset$)
 , $v(\epsilon)=v_{l(\epsilon)}(\epsilon)$  et
$\alpha_{i}(\epsilon)=v_{i}(\epsilon)\mu_{i}$. 
On définit de m\^eme $\displaystyle{\underline{v}(\epsilon) = 
\prod_{ {\tiny \begin{array}{ll}  1\leq k \leq N, \\ 
 k \in \pi_{ +}(\epsilon) \end{array}} } 
\underline{r_{\mu_{k}}} \in \underline{W}}$. 
 On d{\'e}finit un ordre sur $\mathcal{E}$ par : 

$$\epsilon \leq \epsilon' <=> \pi_{+}(\epsilon) \subset
\pi_{+}(\epsilon').$$

On d{\'e}montre alors facilement la proposition suivante : 

\begin{prop}
    
(i) Pour tout $\epsilon$ de $\mathcal{E}$, $Y_{\epsilon}$ est un
espace affine de dimension r{\'e}elle $2l(\epsilon)$. 

(ii) Pour tout $\epsilon \in \mathcal{E}$, $\overline{Y_{\epsilon}} = 
\coprod_{\epsilon' \leq \epsilon} Y_{\epsilon'}$

(iii) $\Gamma = \coprod_{\epsilon \in
  \mathcal{E}} Y_{\epsilon}$

(iv) Pour tout $\epsilon \in \mathcal{E}$, $Y_{\epsilon}$ est stable
par l'action de $T$.
   
\end{prop}

De plus, nous allons avoir besoin du lemme suivant :

\begin{lem}
    
(i) L'ensemble $\Gamma^T$ des points fixes de $\Gamma$ sous l'action de
$T$ est constitu{\'e} des $2^N$ points :

$$[k_{1}, k_{2}, \ldots, k_{N}], \hspace{0,15 cm} o\grave{u} \hspace{0,15 cm}
k_{i} \in \{ e, k_{\mu_{i}} \}.$$

On identifiera donc $\Gamma^T$ avec $\mathcal{E}$ en identifiant $e$
avec $0$ et $k_{\mu_{i}}$ avec $1$.

(ii) Soit $(\epsilon,\epsilon') \in \mathcal{E}^2$, alors : 

$$ \epsilon \in \overline{Y_{\epsilon'}} <=> \epsilon \leq \epsilon',$$

et dans ce cas pour tout élément $h  \in \mathfrak{h}$, on a :  

$${ \rm det}(1-e^{h}|T_{\epsilon'}^{\epsilon})=
\prod_{i \in \pi_{+}(\epsilon')}(1-e^{\alpha_{i}(\epsilon)(h)}),$$
o{\`u} $T_{\epsilon'}^{\epsilon}$ désigne l'espace tangent à
$\overline{Y_{\epsilon'}}$ en
$\epsilon$.

\end{lem}

Comme la variété $\Gamma$ est lisse, $K_{T}(\Gamma)$ s'identifie à $K_{0}(H,
\Gamma)$ (respectivement $K^{0}(H, \Gamma)$) le groupe construit à partir du
semi-groupe des classes d'isomorphisme de faisceaux $H$-équivariants cohérents
(respectivement $H$-équivariants localement libres) sur $\Gamma$. Dans la suite, on
identifiera ces trois groupes. Soit $\epsilon \in \mathcal{E}$ et soit
$\mathcal{F}$ un faisceau $H$-équivariant localement libre sur $\Gamma$, on définit son
caractère $\chi(\overline{Y_{\epsilon}}, \mathcal{F}) \in X[T]$ sur
  $\overline{Y_{\epsilon}}$ par :  
$$\forall t \in T, \hspace{0,1 cm} \chi(\overline{Y_{\epsilon}},
\mathcal{F})(t)=\sum_{k}(-1)^k{ \rm Tr }
 (t; { \rm  H}^k(\overline{Y_{\epsilon}},\mathcal{F}_{/\overline{Y_{\epsilon}}})).$$

La d{\'e}composition $\Gamma=\coprod_{\epsilon \in
  \mathcal{E}}Y_{\epsilon}$ munit $\Gamma$ d'une structure de
  $CW$-complexe $T$-{\'e}quivariant o{\`u} toutes les cellules sont lisses et de
  dimension paire; de plus l'ensemble des points fixes de l'action de $T$ sur
  $\Gamma$ est discret. Grace à cette structure, on a la proposition suivante : 

\begin{prop}

(i) La  $K$-théorie $T$-{\'e}quivariante de $\Gamma^T$ s'identifie {\`a}
l'alg{\`e}bre des fonctions sur $\mathcal{E}$ à valeurs dans $R[T]$, 
qu'on notera $F(\mathcal{E};R[T])$.

(ii) La restriction aux points fixes $i_{T}^*$ : $K_{T}(\Gamma)
 \rightarrow F(\mathcal{E};R[T])$ est injective.
    
(iii) La $K$-théorie $T$-{\'e}quivariante de $\Gamma$ est un $R[T]$-module libre
qui admet comme base la famille $\{\hat{\mu}_{\epsilon}\}_{\epsilon \in
\mathcal{E}}$ caract{\'e}ris{\'e}e par : 

$$\chi(\overline{Y_{\epsilon'}}, \hat{\mu}_{\epsilon})=\delta_{\epsilon', \epsilon}.$$
   
\end{prop}

\begin{defi}
    
Pour $\epsilon \in \mathcal{E}$, on d{\'e}finit $\mu_{\epsilon} \in
F(\mathcal{E};R[T])$
par :

$$\left\{ \begin{array}{ll} \mu_{\epsilon}(\epsilon') = 
\displaystyle{\prod_{i \in \pi_{+}(\epsilon')}e^{\alpha_{i}(\epsilon')}
\prod_{i \in \pi_{+}(\epsilon)}(e^{-\alpha_{i}(\epsilon')}-1) }
 & { \rm si } \hspace{0,15 cm} \epsilon \leq \epsilon', \\ 
 \mu_{\epsilon}(\epsilon') = 0
 & { \rm  sinon }.
\end{array}\right.$$

\end{defi}

On a alors le théorème suivant :

\begin{theo}
    
Pour tout $\epsilon \in \mathcal{E}$, on a : 

$$i_{T}^*(\hat{\mu}_{\epsilon})=\mu_{\epsilon}.$$
   
\end{theo}

\sc 
\begin{flushleft}
D{\'e}monstration :
\end{flushleft} 
\rm

En utilisant la formule de localisation d'Atiyah-Bott (voir \cite{ab}) et le lemme $3$, on
obtient pout tout  $\hat{\mu} \in K_{T}(\Gamma)$ et tout $\epsilon \in \mathcal{E}$: 

\begin{equation} \label{ab} \chi(\overline{Y_{\epsilon}}, \hat{\mu})=  \sum_{\epsilon' \leq
\epsilon}\frac{i_{T}^*(\hat{\mu})(\epsilon')}{\prod_{i \in \pi_{+}(\epsilon)}
(1-e^{\alpha_{i}(\epsilon')})}. \hspace{2 cm} \end{equation}

Soit $\epsilon_{0} \in \mathcal{E}$, et soit
$\mu'_{\epsilon_{0}}=i_{T}^*(\hat{\mu}_{\epsilon_{0}})$. Montrons par
r{\'e}currence sur la longueur de $\epsilon$ que pour tout $\epsilon \in
\mathcal{E}$,
$\mu'_{\epsilon_{0}}(\epsilon)=\mu_{\epsilon_{0}}(\epsilon)$.
Grace {\`a} la
formule \ref{ab} et {\`a} la caract{\'e}risation de $\hat{\mu}_{\epsilon_{0}}$,
on d{\'e}montre facilement par r{\'e}currence sur $l(\epsilon)$ que si 
$\epsilon$ n'est pas plus grand que $\epsilon_{0}$, on
a bien $\mu'_{\epsilon_{0}}(\epsilon)=0$. On peut donc se limiter
au cas o{\`u} $\epsilon_{0} \leq \epsilon$. Si $\epsilon=\epsilon_{0}$, la
formule \ref{ab} et le fait que
$\chi(\overline{Y_{\epsilon_{0}}}, \hat{\mu}_{\epsilon_{0}})=1$ nous
donne bien
$\mu'_{\epsilon_{0}}(\epsilon_{0})=\mu_{\epsilon_{0}}(\epsilon_{0})$.
  Soit $\epsilon > \epsilon_{0}$. On suppose
le r{\'e}sultat v{\'e}rifi{\'e} pour tout $\epsilon'$ de longueur strictement plus
petite que $\epsilon$, on applique la formule \ref{ab} et le fait
que $\chi(\overline{Y_{\epsilon}}, \hat{\mu}_{\epsilon_{0}})=0$  pour
obtenir : 

$$\sum_{\epsilon_{0} \leq  \epsilon' <
\epsilon}\frac{\displaystyle{\prod_{i \in \pi_{+}(\epsilon')}e^{\alpha_{i}(\epsilon')}
\prod_{i \in \pi_{+}(\epsilon_{0})}(e^{-\alpha_{i}(\epsilon')}-1)}}{\displaystyle{\prod_{i \in
  \pi_{+}(\epsilon)}(1-e^{\alpha_{i}(\epsilon')})}}+\frac{\mu'_{\epsilon_{0}}(\epsilon)}
{\displaystyle{\prod_{i \in \pi_{+}(\epsilon)}(1-e^{\alpha_{i}(\epsilon)})} } = 0,$$

d'o{\`u} : 

$$\frac{\mu'_{\epsilon_{0}}(\epsilon)}
{\displaystyle{\prod_{i \in \pi_{+}(\epsilon)}(1-e^{\alpha_{i}(\epsilon)})} }
=-\sum_{\epsilon_{0} \leq  \epsilon' <
\epsilon}\frac{\displaystyle{\prod_{i \in \pi_{+}(\epsilon') \setminus \pi_{+}(\epsilon_{0})}
e^{\alpha_{i}(\epsilon')}} }
{\displaystyle{\prod_{i \in \pi_{+}(\epsilon)\setminus \pi_{+}(\epsilon_{0})}
(1-e^{\alpha_{i}(\epsilon')})}}.$$

Si on pose $\tilde{\epsilon}=\epsilon-(j)$, o{\`u} $j$ est le plus grand
{\'e}l{\'e}ment de $\pi_{+}(\epsilon) \setminus  \pi_{+}(\epsilon_{0})$, on a alors : 

$$\mbox{ \small { \mbox{ $ \frac{ \displaystyle{ \mu'_{\epsilon_{0}}(\epsilon)
  }}{\displaystyle{\prod_{i \in \pi_{+}(\epsilon)}
(1-e^{\alpha_{i}(\epsilon)})}}$}}}\!\!\!\!\!\mbox{ \small { \mbox{ $=\! \frac{ \displaystyle{ \prod_{i \in \pi_{+}(\epsilon) \setminus
  \pi_{+}(\epsilon_{0})} e^{\alpha_{i}(\epsilon) }}}{\displaystyle{\prod_{i \in \pi_{+}(\epsilon)\setminus \pi_{+}(\epsilon_{0})}
(1-e^{\alpha_{i}(\epsilon)})}}-$}}}\!\!\!\!\!\sum_{
 \tiny \begin{array}{ll}  \epsilon_{0} \leq \epsilon' <
\epsilon \\
\hspace{0,25 cm} \epsilon' \neq \tilde{\epsilon} 
\end{array}}\!\!\!\!\!\!\!\!\mbox{ \small { \mbox{ $ \frac{ \displaystyle{ \prod_{i \in \pi_{+}(\epsilon') \setminus
  \pi_{+}(\epsilon_{0})} e^{\alpha_{i}(\epsilon') }}}{\displaystyle{\prod_{i \in \pi_{+}(\epsilon)\setminus \pi_{+}(\epsilon_{0})}
(1-e^{\alpha_{i}(\epsilon')})}}$}}}.
  $$

En effet, comme $j$ est le plus grand {\'e}l{\'e}ment de $\pi_{+}(\epsilon) \setminus
\pi_{+}(\epsilon_{0})$, pour tout $i \in \pi_{+}(\epsilon) \setminus
\pi_{+}(\epsilon_{0})$,  
$\alpha_{i}(\epsilon)=\alpha_{i}(\tilde{\epsilon})$ si $i \neq j$ et 
$\alpha_{j}(\epsilon)=-\alpha_{j}(\tilde{\epsilon})$. On utilise alors la relation 
$\frac{e^{-x}}{1-e^{-x}}=-\frac{e^{x}}{1-e^{x}}$
. Cette m\^eme relation montre, en distinguant les termes qui ont un $1$ en $j${\`e}me position
et ceux qui ont un $0$ en $j${\`e}me position, que la deuxième somme est
nulle et on obtient alors bien : 
$$\mu'_{\epsilon_{0}}(\epsilon)=\displaystyle{\prod_{i \in \pi_{+}(\epsilon)}e^{\alpha_{i}(\epsilon)}
\prod_{i \in \pi_{+}(\epsilon_{0})}(e^{-\alpha_{i}(\epsilon)}-1)
}=\mu_{\epsilon_{0}}(\epsilon).$$

\section{Une autre démonstration du théorème $1$ dans le cas fini}

Si on note $S$ l'algèbre symétrique de $\mathfrak{h}^*$, un résultat analogue à
la proposition $1$ prouvé par A. Arabia dans \cite{aa} montre que la cohomologie
équivariante de $X$ s'identifie à une sous-algèbre $\Lambda$ de $F(W,S)$,
l'algèbre des fonctions de $W$ à valeurs dans $S$ munie de l'addition et de la
multiplication point par point. Cette algèbre $\Lambda$ est un $S$-module libre
qui admet une base $\xi^{w}$ indexée par $W$. Ces fonctions vérifient
$\xi^{w}(v)=0$ si $v$ n'est pas plus grand que $w$ et la formule suivante : 

\begin{theo} Soit $w \leq v$ et soit $v=r_{i_{1}} \cdots r_{i_{l}}$ une
  décomposition réduite de $v$. Si $l(w)=k$, alors : 
$$ \xi^{w}(v) = \sum_{i_{j_{1}} < \cdots < i_{j_{k}}}\beta_{j_{1}}
\cdots \beta_{j_{k}},$$
où la somme porte sur l'ensemble des indices $1 \leq j_{1} < \cdots <
j_{k}\leq l$ tels
que $r_{i_{j_{1}}} \cdots r_{i_{j_{k}}}=w$.

\end{theo}

Ce théorème est démontré de manière purement combinatoire par S. Billey dans
\cite{sb}. Dans la note \cite{mw}, on a retrouvé ce résultat en utilisant les
variétés de Bott-Samelson. Le but de la suite est d'expliquer plus
géométriquement le théorème $1$ à l'aide de ces variétés dans le cas fini.

 Soit $w_{0}=r_{\mu_{1}} \cdots r_{\mu_{N}}$ une d{\'e}composition r{\'e}duite 
du plus grand élément $w_{0}$ de $W$. On pose $\Gamma=
\Gamma(\mu_{1}, \ldots, \mu_{N})$ et on définit une application $T$-équivariante 
$g$ de $\Gamma$ dans $X$ par multiplication (i.e. $g([k_{1},\ldots
,k_{N}])=k_{1}* \cdots *k_{N}
\hspace{0.3cm}  [T]$).

La différence fondamentale avec
ce qui se passe en cohomologie est due au fait que en cohomologie les 
opérateurs de Demazure $A_{i}$ vérifient $A_{i}^2=0$ alors qu'en $K$-théorie, les
opérateurs $D_{i}$ vérifient $D_{i}^2=D_{i}$. Une généralisation immédiate de la
proposition $3.36$ de \cite{kkk}, où le résultat n'est énoncé que pour des
décompositions réduites, nous donne : 
$$\forall \tau \in K_{T}(X) \hspace{0,1 cm}, \hspace{0,2 cm} \chi(
\overline{Y}_{\epsilon},g^*(*\tau))=*(D_{\underline{v}(\epsilon)}i_{T}^*(\tau))(1),$$
où pour $\underline{v} \in \underline{W}$, on a posé $D_{\underline{v}}=D_{v}$. Or d'après les
propriétés $(i)$, $(ii)$ et $(iii)$ de la proposition $1$ : 
$$\forall (v,w) \in W^2, (D_{v}(\psi^{w}))(1)=\delta_{v,w}.$$

On déduit des deux formules précédentes que pour tout $\epsilon \in \mathcal{E}$ et
tout élément $w \in W$, on a
\begin{equation} \label{chi} \chi(\overline{Y}_{\epsilon},g^*(*\hat{\psi}^w))=\delta_{\underline{v}(\epsilon)
 ,\underline{w}}. \hspace{2 cm} \end{equation}

D'après la caractérisation de la base
$\{\hat{\mu}_{\epsilon}\}_{\epsilon \in \mathcal{E}}$, on a donc :

$$\forall w \in W, g^*(\hat{\psi}^{w})=\sum_{\epsilon \in \mathcal{E},
  \underline{v}(\epsilon)=\underline{w}}*\hat{\mu}_{\epsilon} \hspace{0,1 cm}.$$
 
Soit $(w,v) \in W^2$. Si on choisit un élément $\epsilon' \in \mathcal{E}$ tel que
$g(\epsilon')=v$ et tel que $l(\epsilon')=l(v)$ (ce qui correspond au choix d'une
décomposition réduite de $v$), la formule précédente nous montre que
$$\displaystyle{\psi^{w}(v)=(i_{T}^*(
\hat{\psi}^{w}))(v)=\sum_{\epsilon \in \mathcal{E},
  \underline{v}(\epsilon)=\underline{w}}*\mu_{\epsilon}(\epsilon')},$$
 ce qui nous redonne
bien le théorème $1$ à l'aide du théorème $3$. On remarque que l'indépendance
par rapport au choix d'une décomposition réduite de $v$ est ici une conséquence
immédiate du fait que pour tout fibré $T$-équivariant $\tau$ sur $X$, la représentation
de $T$ dans $g^*(\tau)_{x}$ est
la m\^eme en tout point $x$ de $\Gamma$ qui vérifie $g(x)=v$. Cette
vision géométrique montre également que la formule du théorème $1$ reste
valable si on prend une décomposition non réduite de $v$ (ce qu'on peut aussi
voir de manière combinatoire).

\section{Une autre base de $K_{T}(X)$}

Comme $X$ est une variété lisse, tout comme dans le cas des variétés de
Bott-Samelson, $K_{T}(X)$ s'identifie à $K_{0}(H,X)$ et à $K^0(H,X)$. On sait
alors que la décomposition en cellules de Schubert $X=\coprod_{w \in
  W}X_{w}$ fournit une base $\{[\mathcal{O}]_{\overline{X}_{w}}\}_{w \in W}$ de
$K_{0}(X,H)$. Les classes $[\mathcal{O}]_{\overline{X}_{w}}$ sont définies par le faisceau
structural de $\overline{X}_{w}$ prolongé par $0$ sur $X \setminus
\overline{X}_{w}$. Pour $w \in W$, on pose
$\hat{\sigma}^w=*[\mathcal{O}]_{\overline{X}_{w}} \in K_{T}(X)$, et
$\sigma^{w}=i_{T}^*(\hat{\sigma}^w)$.
 Le résultat suivant est prouvé dans
\cite{kkk} : 

\begin{prop}

Pour tout $w \in W$ et tout entier $1 \leq i \leq r$,

$$D_{i}(\sigma^w)=\left\{ \begin{array}{ll} \sigma^w 

 & {\rm si } \hspace{0.2 cm} wr_{i}<w, \\ 
\sigma^{wr_{i}}
 & {\rm si } \hspace{0.2 cm} wr_{i}>w.  \end{array} \right.$$

\end{prop}

On définit les éléments $a_{w}^v \in R[T]$ par 
$\displaystyle{\sigma^{w}=\sum_{v \in W}a_{w}^v\psi^v}$. On va donner une
expression explicite de ces coefficients. Soit $w \in W$ et soit $r_{i}$ 
une réflexion simple telle que $wr_{i}>w$. Si on
applique l'opérateur $D_{i}$ à la décomposition 
$\displaystyle{\sigma^{w}=\sum_{v \in W}a_{w}^v\psi^v}$, on
obtient : 
$$\sigma^{wr_{i}}=D_{i}(\sigma^w)=
\sum_{v \in W}a_{w}^v D_{i}\psi^{v}.$$

En utilisant les relations vérifiées par les fonctions $\psi^{v}$, on trouve
alors : 
$$\sum_{v \in W}a_{wr_{i}}^{v}\psi^{v}=\sum_{v \in W, vr_{i}<v}a_{w}^{v}(\psi^{v} + 
\psi^{vr_{i}})=\sum_{v \in W, vr_{i}<v}a_{w}^{v}\psi^{v}
+ \sum_{v \in W, vr_{i}>v}a_{w}^{vr_{i}}\psi^{v}. $$

On obtient donc la relation de récurrence suivante sur les coefficients
$a_{w}^{v}$ : 

$$a_{wr_{i}}^{v}=\left\{ \begin{array}{ll} a_{w}^{v}
 & {\rm si } \hspace{0.2 cm} vr_{i}<v, \\ 
a_{w}^{vr_{i}}
 & {\rm si } \hspace{0.2 cm} vr_{i}>v.  \end{array} \right.$$

De ces relations, on déduit en utilisant les relations \ref{hecke1} et
\ref{hecke2} : 

\begin{equation} \label{rec} 
\forall (w,v) \in W^2, \quad
a_{w}^{v}=a_{1}^{\underline{v} \hspace{0,1 cm} \underline{w^{-1}}} ,
 \end{equation}

 où pour $\underline{v} \in \underline{W}$, on a posé 
$a_{1}^{\underline{v}}=a_{1}^{v}$. Il suffit donc de trouver la décomposition de
$\sigma^{1}$. Pour cela, on aura besoin des valeurs de $\sigma^{1}$ : 

$$\left\{ \begin{array}{l} \sigma^{1}(1)=\prod_{\alpha \in
      \Delta_{+}}
(1-e^{-\alpha}),
   \\ 
\sigma^{1}(v)=0
 \hspace{0,5 cm} {\rm si } \hspace{0.2 cm} v \neq 1.  \end{array} \right.$$
La valeur de $ \sigma^{1}(1)$ est calculée à l'aide de la formule
d'auto-intersection et les autres valeurs sont nulles par le théorème de localisation.

Comme on a $\hat{\sigma}^{1}=\sum_{v \in
    W}a_{1}^{v}\hat{\psi}^{v}$, soit $v \in W$ et $\epsilon \in \mathcal{E}$ 
tel que $u(\epsilon)=v$, d'après la formule \ref{chi}, le coefficient $a_{1}^{v}$
est donné par :
$$a_{1}^{v}=*\chi(\overline{Y}_{\epsilon},g^*(*\hat{\sigma}^{1})).$$

En utilisant la formule \cite{ab} et les valeurs de
$*\sigma^{1}=*i_{T}^*\hat{\sigma}^{1}$, on obtient alors : 
$$a_{1}^{v}=  \sum_{\epsilon' \leq
\epsilon, u(\epsilon')=1}\frac{\prod_{\alpha \in
      \Delta_{+}}
(1-e^{-\alpha})}{\prod_{i \in \pi_{+}(\epsilon)}
(1-e^{-\alpha_{i}(\epsilon')})}. $$

On a donc la proposition suivante : 

\begin{prop}

Soit $v \in W$ et soit $v=r_{i_{1}} \cdots r_{i_{l}}$ une décomposition réduite de
$v$. Pour tout sous-ensemble $I$ de $\{1, \ldots , l \}$ et tout entier $1 \leq i
\leq l$, on pose $\displaystyle{\beta_{i}(I)=(\prod_{j \in I , j \leq
  i}r_{i_{j}})\alpha_{i}}$ ($\beta_{i}(I)=\alpha_{i}$ si $I \cap \{ 1, \ldots ,
i \} = \emptyset$). Alors 
$$\displaystyle{\sum_{k=0}^{l}\sum
  \frac{\prod_{\alpha \in
  \Delta_{+}}(1-e^{-\alpha})
}{
\prod_{i=1}^{l}
(1-e^{-\beta_{i}(\{j_{1},\ldots , j_{k}\}  )})}},$$
 où la deuxième somme porte sur l'ensemble des indices 
$1 \leq j_{1} < \cdots < j_{k}\leq l$ tels que $r_{i_{j_{1}}} \cdots r_{i_{j_{k}}}=1$, est
un élément de $R[T]$ qui ne dépend pas du choix d'une décomposition réduite de
$v$. Si on note $b^{v}$ cet élément, alors $a_{1}^{v}=b^{v}$.
 
\end{prop}

\begin{Rmq}

La proposition précédente est encore valable si on prend une décomposition non
réduite de $v$.

\end{Rmq}

Si on utilise la relation \ref{rec}, et si pour $\underline{v} \in \underline{W}$, on pose
$b^{\underline{v}}=b^v$, on obtient alors le théorème suivant

\begin{theo} 

Soit $w \in W$, alors : 
$$*[\mathcal{O}_{\overline{X}_{w}}]=\sum_{v \in W}
b^{\underline{v}\hspace{0,1 cm} \underline{w^{-1}} }\hat{\psi}^{v}.$$

\end{theo}

Soit $Q_{W}$ le $Q[T]$-module libre qui admet pour base la famille $\{\delta_{w}\}_{w
  \in W}$ et qu'on munit d'une structure d'anneau définie par : 
$$(q_{1}\delta_{w_{1}}).(q_{2}\delta_{w_{2}})=q_{1}(w_{1}q_{2})\delta_{w_{1}w_{2}}
\hspace{0,1 cm}, \hspace{0,1 cm} \forall (q_{1}, q_{2}) \in Q[T]^2 , \hspace{0,1 cm} { \rm et }
\hspace{0,1 cm}  (w_{1},w_{2}) \in W^2,$$
où l'action de $W$ sur $Q[T]$ est déduite de celle de $W$ sur $T$. Dans
\cite{kkk}, B. Kostant et S. Kumar introduisent des éléments $\{y_{i} \}_{1 \leq
  i \leq r}$ de $Q[T]$ définis par : 
$$y_{i}=\frac{1}{1-e^{-\alpha_{i}}}(\delta_{1}-e^{-\alpha_{i}}\delta_{r_{i}}).$$

Les $y_{i}$ vérifiant les relations de tresses, on peut
définir un élément $y_{w} \in Q_{W}$ pour tout $w \in W$. 
On définit alors des éléments $\{ b_{v,w}
\}_{(v,w) \in W^2}$ de $Q[T]$ par  : 
$$y_{v^{-1}}=\sum_{w \in W}b_{v,w}\delta_{w^{-1}}.$$
 
D'après l'expression combinatoire de $b_{v,w}$ donnée par le lemme $3.5$ de \cite{ku},  
$\displaystyle{b^v=b_{v^{-1},1}\prod_{\alpha \in \Delta_{+}}(1-e^{-\alpha})}$, 
et de plus dans \cite{ku}, S. Kumar montre que quand  $\overline{X}_{v}$
est lisse, $\displaystyle{b_{v,1}=\prod_{\gamma \in S(v)}(1-e^{-\gamma})^{-1}}$, 
où pour $u \in W$, $S(u)=\{ \alpha \in R^{+}, r_{\alpha} \leq u \}$. En
particulier $b^{w_{0}}=1$ et donc d'après le théorème $5$ : 
$$*[\mathcal{O}_{\overline{X}_{w_{0}}}]=\sum_{w \in W}\hat{\psi}_{w}.$$

 \bibliography{Ktheory}
 \bibliographystyle{smfplain}

\vspace{1 cm}

\bf{e-mail}  \rm : willems@math.jussieu.fr

\end{document}